\documentclass[reqno]{article}%
\usepackage{amsmath}
\usepackage{graphicx}
\usepackage{amsfonts}
\usepackage{amssymb}%
\newtheorem{theorem}{Theorem}

\newtheorem{conjecture}[theorem]{Conjecture}

\newtheorem{lemma}[theorem]{Lemma}

\newtheorem{proposition}[theorem]{Proposition}

\newenvironment{proof}[1][Proof]{\noindent{\textbf {#1}  }}  {\hfill$\Box$\bigskip}

\begin{document}

\title{Graphs and Hermitian matrices: discrepancy and singular values}
\author{B\'{e}la Bollob\'{a}s\thanks{Department of Mathematical Sciences, University
of Memphis, Memphis TN 38152, USA} \thanks{Trinity College, Cambridge CB2 1TQ,
UK} \thanks{Research supported by NSF grant ITR 0225610 and DARPA grant
F33615-01-C-1900 } \ and Vladimir Nikiforov$^{\ast}$\thanks{Contact author.
E-mail address: \textit{vnikifrv@memphis.edu}}}
\maketitle

\begin{abstract}
Let $A=\left(  a_{ij}\right)  _{i,j=1}^{n}$ be a Hermitian matrix of size
$n\geq2$, and set%
\begin{align*}
\rho\left(  A\right)   &  =\frac{1}{n^{2}}\sum_{i=1}^{n}\sum_{j=1}^{n}%
a_{ij},\\
disc\left(  A\right)   &  =\max_{X,Y\subset\left[  n\right]  ,X\neq
\varnothing,Y\neq\varnothing}\frac{1}{\sqrt{\left\vert X\right\vert \left\vert
Y\right\vert }}\left\vert \sum_{i\in X}\sum_{j\in Y}\left(  a_{ij}-\rho\left(
A\right)  \right)  \right\vert .
\end{align*}
We show that the second singular value $\sigma_{2}\left(  A\right)  $ of $A$
satisfies%
\[
\sigma_{2}\left(  A\right)  \leq C_{1}disc\left(  A\right)  \log n,
\]
for some absolute constant $C_{1},$ and this is best possible up to a
multiplicative constant. Moreover, we construct infinitely many dense regular
graphs $G$ such that
\[
\sigma_{2}\left(  A\left(  G\right)  \right)  \geq C_{2}disc\left(  A\left(
G\right)  \right)  \log\left\vert G\right\vert
\]
where $C_{2}>0$ is an absolute constant and $A\left(  G\right)  $ is the
adjacency matrix of $G.$ In particular, these graphs disprove two conjectures
of Fan Chung.

\textbf{Keywords: }\textit{discrepancy, graph eigenvalues, second singular
value, pseudo-random graphs, quasi-random graphs}

\end{abstract}

\section{Introduction}

Given a Hermitian matrix $A$ of size $n,$ let $\mu_{1}\left(  A\right)
\geq...\geq\mu_{n}\left(  A\right)  $ be its eigenvalues, and $\sigma
_{1}\left(  A\right)  \geq...\geq\sigma_{n}\left(  A\right)  $ be its singular
values. As $A$ is Hermitian, the values $\sigma_{i}\left(  A\right)  $ are the
moduli of $\mu_{i}\left(  A\right)  $ taken in descending order, so
$\sigma_{1}\left(  A\right)  $ is the $l_{2}$ operator norm and the spectral
radius of $A$.

Our graph-theoretic notation is standard (e.g., see \cite{Bol1}). For
simplicity, graphs are assumed to be defined on the vertex set $\left[
n\right]  =\left\{  1,...,n\right\}  .$ Occasionally, to remind the reader of
this, we write $G\left(  n\right)  $ for a graph of order $n,$ $G\left(
n,m\right)  $ for a graph of order $n$ and size $m$. We write $e\left(
X\right)  $ for $e\left(  G\left[  X\right]  \right)  $ if it is understood
which graph $G$ is to be taken. Given a graph $G,$ we let $\mu_{i}\left(
G\right)  =\mu_{i}\left(  A\left(  G\right)  \right)  ,$ and $\sigma
_{i}\left(  G\right)  =\sigma_{i}\left(  A\left(  G\right)  \right)  ,$ where
$A\left(  G\right)  $ is the adjacency matrix of $G.$ Given a graph $G$ and
$X,Y\subset V\left(  G\right)  ,$ we denote by $e\left(  X,Y\right)  $ the
number of the ordered pairs $\left(  u,v\right)  $ such that $u\in X,$ $v\in
Y,$ and $u$ is adjacent to $v$.

For every graph $G=G\left(  n\right)  ,$ set $\rho\left(  G\right)  =e\left(
G\right)  \binom{n}{2}^{-1}$ and let
\begin{align*}
disc_{1}\left(  G\right)   &  =\max_{X\subset V\left(  G\right)
,X\neq\varnothing}\left\{  \frac{1}{\left\vert X\right\vert }\left\vert
e\left(  X\right)  -\rho\left(  G\right)  \binom{\left\vert X\right\vert }%
{2}\right\vert \right\}  ,\\
disc_{2}\left(  G\right)   &  =\max_{X,Y\subset V\left(  G\right)
,X\neq\varnothing,Y\neq\varnothing}\left\{  \frac{1}{\sqrt{\left\vert
X\right\vert \left\vert Y\right\vert }}\left\vert e\left(  X,Y\right)
-\rho\left(  G\right)  \left\vert X\right\vert \left\vert Y\right\vert
\right\vert \right\}  .
\end{align*}

The function $disc_{1}\left(  G\right)  $ is, in fact, Thomason's coefficient
$\alpha$ in his definition of $\left(  p,\alpha\right)  $-jumbled graphs (see,
e.g. \cite{Tho1}, \cite{Tho2}), on which he based his study of pseudo-random
graphs.\ In principle, $disc_{2}\left(  G\right)  $ has the same role as
$disc_{1}\left(  G\right)  ,$ although the two invariants may differ
significantly for certain graphs, e.g., the star $K_{1,n}$. Chung, Graham, and
Wilson \cite{CGW} (see also \cite{KrSu}) used coarser functions to describe
the edge distribution of a graph, thus introducing the quasi-random graph
properties. Surprisingly, these properties can be expressed in terms of the
two largest moduli of the eigenvalues (or equivalently, the two largest
singular values) of the adjacency matrix of a graph; we refer the interested
reader to \cite{KrSu} for more details. A natural question is, whether similar
relations exist between singular values and the functions $disc_{1}\left(
G\right)  $ and $disc_{2}\left(  G\right)  .$

Chung (\cite{Ch1}, p. 35) made the following interesting conjecture concerning
$\sigma_{2}\left(  G\right)  $ and $disc_{2}\left(  G\right)  .$

\begin{conjecture}
\label{con1} There is an absolute constant $C$ such that for every regular
graph $G,$%
\begin{equation}
\sigma_{2}\left(  G\right)  <Cdisc_{2}\left(  G\right)  . \label{d1}%
\end{equation}

\end{conjecture}

The main goal of this paper is to study similar questions for graphs and
Hermitian matrices. In particular, in section \ref{sec2} we define a function
$disc\left(  A\right)  $ for a Hermitian matrix $A$ that naturally extends the
function $disc_{2}\left(  G\right)  ,$ and show that there is some constant
$C^{\prime}$ such that for every Hermitian matrix $A$ of size $n\geq2,$ we
have
\begin{equation}
\sigma_{2}\left(  A\right)  <C^{\prime}disc\left(  A\right)  \log n.
\label{upbnd}%
\end{equation}
We explicitly construct a nonnegative symmetric matrix showing that
(\ref{upbnd}) is best possible up to a multiplicative constant. Moreover, in
section \ref{sec3} we construct infinitely many dense regular graphs $G$ such
that
\[
\mu_{2}\left(  G\right)  >C^{\prime\prime}disc_{2}\left(  G\right)
\log\left\vert G\right\vert
\]
for some absolute constant $C^{\prime\prime}>0,$ thus disproving Conjecture
\ref{con1}. In fact, as we show in section \ref{sec3.last}, these graphs
disprove also another conjecture of Chung stating a similar problem for
Laplacian eigenvalues (\cite{Ch2}, p. 77).

In particular, in section \ref{Tho} we show that the bound on $disc_{1}\left(
G\right)  $ due to Thomason (\cite{Tho3}, Theorem 1) can be easily extended to
$disc_{2}\left(  G\right)  .$

Recently Chung and Graham discussed in \cite{CG3} quasi-random graph
properties of sparse graphs. In particular, they asked whether for sparse
graphs small discrepancy implies small second singular value, as in the case
of dense graphs. Krivelevich and Sudakov gave an explicit example in
\cite{KrSu} that answers this question in the negative. In section \ref{sec4}
we describe a general construction showing that such examples are not exceptional.

\section{\label{sec2} Second singular value and discrepancy of Hermitian
matrices}

Given a matrix $A=\left(  a_{ij}\right)  _{i,j=1}^{n}$ and nonempty sets
$I,J\subset\left[  n\right]  ,$ we denote by $A\left[  I,J\right]  $ the
submatrix of the entries $a_{ij}$ with $i\in I,$ $j\in J$. We write $E_{n}$
for the $n\times n$ matrix of all ones, and denote by $\left\langle
\mathbf{x,y}\right\rangle $ the standard inner product of two vectors
$\mathbf{x,y}\in\mathbb{C}^{n}$. For a Hermitian matrix $A=\left(
a_{ij}\right)  _{i,j=1}^{n}$ set
\begin{align}
\rho^{\prime}\left(  A\right)   &  =\frac{1}{n^{2}}\sum_{i=1}^{n}\sum
_{j=1}^{n}a_{ij},\nonumber\\
disc\left(  A\right)   &  =\max_{X,Y\subset\left[  n\right]  ,X\neq
\varnothing,Y\neq\varnothing}\frac{1}{\sqrt{\left|  X\right|  \left|
Y\right|  }}\left|  \sum_{i\in X}\sum_{j\in Y}\left(  a_{ij}-\rho^{\prime
}\left(  A\right)  \right)  \right|  . \label{defdisc}%
\end{align}

In this section we investigate the relationship between $\sigma_{2}\left(
A\right)  $ and $disc\left(  A\right)  .$ Our main goal is to prove Theorem
\ref{th1} and to show that the assertion is best possible up to a
multiplicative constant.

Observe that for a graph $G$ the value $\rho\left(  G\right)  $ is, generally
speaking, different from $\rho^{\prime}\left(  A\left(  G\right)  \right)  $
of the adjacency matrix $A\left(  G\right)  $ of $G.$ However, setting
$A=A\left(  G\right)  $, we easily see that for every two nonempty sets
$X,Y\subset V\left(  G\right)  ,$%
\begin{align*}
&  \frac{1}{\sqrt{\left\vert X\right\vert \left\vert Y\right\vert }}\left\vert
\left\vert \sum_{i\in X}\sum_{j\in Y}\left(  a_{ij}-\rho^{\prime}\left(
A\right)  \right)  \right\vert -\left\vert \sum_{i\in X}\sum_{j\in Y}\left(
a_{ij}-\rho\left(  G\right)  \right)  \right\vert \right\vert \\
&  \leq\frac{1}{\sqrt{\left\vert X\right\vert \left\vert Y\right\vert }%
}\left\vert \rho^{\prime}\left(  A\right)  -\rho\left(  G\right)  \right\vert
\left\vert X\right\vert \left\vert Y\right\vert \\
&  =\sqrt{\left\vert X\right\vert \left\vert Y\right\vert }\left\vert \frac
{1}{n^{2}}-\frac{1}{n\left(  n-1\right)  }\right\vert 2e\left(  G\right)
\leq\frac{2e\left(  G\right)  }{n\left(  n-1\right)  }.
\end{align*}
Therefore,
\begin{equation}
\left\vert disc_{2}\left(  G\right)  -disc\left(  A\right)  \right\vert
\leq\frac{2e\left(  G\right)  }{n\left(  n-1\right)  }\leq1, \label{mag}%
\end{equation}
i.e., the function $disc\left(  A\right)  $ closely approximates
$disc_{2}\left(  G\right)  $.

\subsection{\label{sec2.1} An upper bound on $\sigma_{2}$}

Observe that for every Hermitian matrix $A$ of size $n,$ and every
$\mathbf{x}\in\mathbb{C}^{n},$ the Rayleigh principle states that
\[
\mu_{n}\left(  A\right)  \left\|  \mathbf{x}\right\|  ^{2}\leq\left\langle
A\mathbf{x,x}\right\rangle \leq\mu_{1}\left(  A\right)  \left\|
\mathbf{x}\right\|  ^{2};
\]
thus, we see that for every $\mathbf{x}\in\mathbb{C}^{n},$
\begin{equation}
\left|  \left\langle A\mathbf{x,x}\right\rangle \right|  \leq\left\|
A\right\|  \left\|  \mathbf{x}\right\|  =\sigma_{1}\left(  A\right)  \left\|
\mathbf{x}\right\|  ^{2}. \label{basin}%
\end{equation}

\begin{theorem}
\label{th1} There is some constant $C$ such that for every Hermitian matrix
$A$ of size $n\geq2,$%
\[
\sigma_{2}\left(  A\right)  \leq Cdisc\left(  A\right)  \log n,
\]

\end{theorem}

Before proceeding to the proof of this theorem let us prove some technical
results. We shall prove first a curious lemma that is somewhat stronger than needed.

\begin{lemma}
\label{lapp} Let $p\geq1$, $n\geq1$ and $0<\varepsilon<1$. Then for every
$\mathbf{x}=(x_{i})_{1}^{n}\in\mathbb{C}^{n}$ with $\Vert\mathbf{x}\Vert
_{p}=1$, there is a vector $\mathbf{y}=(y_{i})_{1}^{n}\in\mathbb{C}^{n}$ such
that $y_{i}$ take no more than
\[
\left\lceil \frac{8\pi}{\varepsilon}\right\rceil \left\lceil \frac
{4}{\varepsilon}\log\frac{4n}{\varepsilon}\right\rceil
\]
values and $\Vert\mathbf{x}-\mathbf{y}\Vert_{p}\leq\varepsilon$.
\end{lemma}

\begin{proof}
We shall prove first that if $x_{i}$ are nonnegative reals, then there is a
$\mathbf{y}=(y_{i})_{1}^{n}$ such that $y_{i}$ are nonnegative reals,
$\Vert\mathbf{x}-\mathbf{y}\Vert_{p}\leq\varepsilon,$ and $y_{i}$ take no more
than
\[
k=\left\lceil \frac{2}{\varepsilon}\log\frac{2n}{\varepsilon}\right\rceil
\]
different values. We may and shall assume $x_{1}\geq...\geq x_{n}\geq0.$

Let us define a sequence $n_{1}<n_{2}<\ldots<n_{l}\leq n$ as follows. Set
$n_{1}=1;$ having defined $n_{i}$, let $s$ be the maximal index such that
\[
x_{s}\geq\left(  1-\frac{\varepsilon}{2}\right)  x_{n_{i}};
\]
if $i=k+1$ or $s=n$, stop the sequence; otherwise, let $n_{i+1}=s+1.$ Finally,
for $1\leq j\leq n$, set $y_{j}=x_{n_{i+1}-1}$ if $n_{i}\leq j<n_{i+1}$ and
$y_{j}=0$ if $n_{k}<j\leq n$. For the sake of convenience, set $n_{l+1}=n+1$.
Then, if $l\leq k,$%
\[
\sum\limits_{j=1}^{n}\left|  x_{j}-y_{j}\right|  ^{p}\leq\sum\limits_{h=1}%
^{l}\sum\limits_{j=n_{h}}^{n_{h+1}-1}\left(  \frac{\varepsilon}{2}%
x_{n_{h+1}-1}\right)  ^{p}\leq\left(  \frac{\varepsilon}{2}\right)  ^{p}%
\sum\limits_{j=1}^{n}x_{j}^{p}=\left(  \frac{\varepsilon}{2}\right)
^{p}<\varepsilon^{p}.
\]
Let now $l\geq k+1;$ observe, that the choice of $k$ implies for every
$j=n_{k+1},...,n,$
\[
x_{j}\leq\left(  1-\frac{\varepsilon}{2}\right)  ^{k}x_{1}\leq\left(
1-\frac{\varepsilon}{2}\right)  ^{k}\leq\frac{\varepsilon}{2n}.
\]
Hence, for $l\geq k+1,$ we obtain,%
\begin{align*}
\sum\limits_{j=1}^{n}\left|  x_{j}-y_{j}\right|  ^{p}  &  =\sum\limits_{j=1}%
^{n_{k+1}-1}\left|  x_{j}-y_{j}\right|  ^{p}+\sum\limits_{j=n_{k+1}}^{n}%
x_{j}^{p}\\
&  \leq\sum\limits_{h=1}^{k}\sum\limits_{j=n_{h}}^{n_{h+1}-1}\left(
\frac{\varepsilon}{2}x_{n_{h+1}-1}\right)  ^{p}+n\left(  \frac{\varepsilon
}{2n}\right)  ^{p}\leq\left(  \frac{\varepsilon}{2}\right)  ^{p}+\left(
\frac{\varepsilon}{2}\right)  ^{p}\leq\varepsilon^{p}.
\end{align*}
Consequently, $\Vert\mathbf{x}-\mathbf{y}\Vert_{p}\leq\varepsilon$, as required.

Let now $\mathbf{x}=(x_{j})_{1}^{n}\in\mathbb{C}^{n}$ be an arbitrary vector
and for every $j\in\left[  n\right]  ,$ let
\[
x_{j}=\left|  x_{j}\right|  \exp\left(  \theta_{j}2\pi i\right)  ,
\]
where $0\leq\theta_{j}<1$. Set
\[
\mathbf{x}^{\ast}\mathbf{=}\left(  \left|  x_{1}\right|  ,\ldots,\left|
x_{n}\right|  \right)  .
\]
According to the above, there exists $\mathbf{z}=(z_{j})_{1}^{n},$ such that
$\Vert\mathbf{x}^{\ast}-\mathbf{z}\Vert_{p}\leq\varepsilon/2,$ $z_{j}\geq0$
and $z_{j}$ take at most%
\[
k=\left\lceil \frac{4}{\varepsilon}\log\frac{4n}{\varepsilon}\right\rceil
\]
different values. Let
\[
m=\left\lceil \frac{8\pi}{\varepsilon}\right\rceil .
\]
We shall show that the vector $\mathbf{y=}(y_{i})_{1}^{n}$ defined by
\[
\left|  y_{j}\right|  =z_{j},\text{ }\arg\left(  y_{j}\right)  =\frac
{\left\lfloor m\theta_{j}\right\rfloor }{m}2\pi
\]
is as required. Let us first check that $\Vert\mathbf{x}-\mathbf{y}\Vert
_{p}\leq\varepsilon.$ Indeed, define $\mathbf{z}^{\ast}=(z_{j}^{\ast})_{1}%
^{n}$ by
\[
z_{j}^{\ast}=z_{i}\exp\left(  \theta_{j}2\pi i\right)  .
\]
We have, by the triangle inequality,
\begin{equation}
\Vert\mathbf{x}-\mathbf{y}\Vert_{p}\leq\Vert\mathbf{x}-\mathbf{z}^{\ast}%
\Vert_{p}+\Vert\mathbf{y}-\mathbf{z}^{\ast}\Vert_{p}. \label{in2}%
\end{equation}
On the one hand,%
\begin{equation}
\sum_{j=1}^{n}\left|  x_{j}-z_{j}^{\ast}\right|  ^{p}=\sum_{j=1}^{n}\left|
\left|  x_{j}\right|  -z_{j}\right|  ^{p}\exp\left(  \theta_{j}2\pi i\right)
^{p}=\sum_{j=1}^{n}\left|  \left|  x_{j}\right|  -z_{j}\right|  ^{p}%
\leq\left(  \frac{\varepsilon}{2}\right)  ^{p}. \label{in1}%
\end{equation}
On the other hand,
\[
\left|  \exp\left(  \frac{\left\lfloor m\theta_{j}\right\rfloor }{m}2\pi
i\right)  -\exp\left(  \theta_{j}2\pi i\right)  \right|  \leq2\sin\left(
\frac{1}{2m}2\pi\right)  <\frac{2\pi}{m}\leq\frac{\varepsilon}{4},
\]
and hence,%
\begin{align*}
\sum_{j=1}^{n}\left|  y_{j}-z_{j}^{\ast}\right|  ^{p}  &  \leq\sum_{j=1}%
^{n}\left|  z_{j}\exp\left(  \theta_{j}2\pi i\right)  -z_{j}\arg\left(
\frac{\left\lfloor m\theta_{j}\right\rfloor }{m}2\pi i\right)  \right|
^{p}\leq\sum_{j=1}^{n}\left(  \frac{\varepsilon}{4}\right)  ^{p}\left|
z_{j}\right|  ^{p}\\
&  =\varepsilon^{p}\left(  \left\|  \mathbf{z}\right\|  _{p}\right)  ^{p}%
\leq\left(  \frac{\varepsilon}{4}\right)  ^{p}\left(  \frac{\varepsilon}%
{2}+\left\|  \mathbf{x}\right\|  _{p}\right)  ^{p}<\left(  \frac{\varepsilon
}{2}\right)  ^{p}.
\end{align*}
Hence, in view of (\ref{in2}) and (\ref{in1}), we obtain
\[
\Vert\mathbf{x}-\mathbf{y}\Vert_{p}\leq\left(  \left(  \frac{\varepsilon}%
{2}\right)  ^{p}+\left(  \frac{\varepsilon}{2}\right)  ^{p}\right)  ^{1/p}%
\leq\varepsilon.
\]
To complete the proof, observe that $y_{i}$ take at most
\[
km\leq\left\lceil \frac{8\pi}{\varepsilon}\right\rceil \left\lceil \frac
{4}{\varepsilon}\log\frac{4n}{\varepsilon}\right\rceil
\]
different values.
\end{proof}

We say that a partition $X=\cup_{i=1}^{m}P_{i}$ is \emph{proper} if the sets
$P_{i}$ are nonempty.

\begin{lemma}
\label{le7} Let $B=\left(  b_{ij}\right)  _{i,j=1}^{n}$ be a Hermitian matrix
and $\left[  n\right]  =\cup_{i=1}^{m}P_{i}$ be a proper partition. Let
$\mathbf{y}\in\mathbb{C}^{m},$ and $\mathbf{x=}\left(  x_{i}\right)  _{i}%
^{n}\in\mathbb{C}^{n}$ be such that $x_{i}=y_{j}$ for every $i\in P_{j}$. Then
the Hermitian matrix $C=\left(  c_{ij}\right)  _{i,j=1}^{m}$ defined by
\[
c_{ij}=\frac{1}{\sqrt{\left\vert P_{i}\right\vert \left\vert P_{j}\right\vert
}}\sum_{r\in P_{i}}\sum_{s\in P_{j}}b_{rs}%
\]
satisfies%
\[
\left\vert \left\langle B\mathbf{x,x}\right\rangle \right\vert \leq\sigma
_{1}\left(  C\right)  \left\Vert \mathbf{x}\right\Vert ^{2}.
\]

\end{lemma}

\begin{proof}
For every $k\in\left[  m\right]  ,$ set $t_{k}=\sqrt{\left|  P_{k}\right|
}y_{k},$ and let $\mathbf{t=}\left(  t_{1},...,t_{m}\right)  ,$ so that
$\left\|  \mathbf{t}\right\|  =\left\|  \mathbf{x}\right\|  .$ Also, we see
that
\begin{align*}
\left\langle B\mathbf{x,x}\right\rangle  &  =\sum_{i=1}^{n}\sum_{j=1}%
^{n}b_{ij}x_{i}\overline{x}_{j}=\sum_{i=1}^{m}\sum_{j=1}^{m}\frac
{t_{i}\overline{t}_{j}}{\sqrt{\left|  P_{i}\right|  \left|  P_{j}\right|  }%
}\sum_{r\in P_{i}}\sum_{s\in Pj}b_{rs}\\
&  =\sum_{i=1}^{m}\sum_{j=1}^{m}c_{ij}t_{i}\overline{t}_{j}=\left\langle
C\mathbf{t,t}\right\rangle .
\end{align*}
Hence, from (\ref{basin}), we obtain
\[
\left|  \left\langle B\mathbf{x,x}\right\rangle \right|  \leq\sigma_{1}\left(
C\right)  \left\|  \mathbf{t}\right\|  ^{2}=\sigma_{1}\left(  C\right)
\left\|  \mathbf{x}\right\|  ^{2},
\]
completing the proof.
\end{proof}

\begin{proof}
[Proof of Theorem \ref{th1}]Set $\rho^{\prime}=\rho^{\prime}\left(  A\right)
$ and let
\[
B=A-\rho^{\prime}E_{n}.
\]
Our first goal is to show that
\begin{equation}
\sigma_{2}\left(  A\right)  \leq\sigma_{1}\left(  B\right)  . \label{eq1}%
\end{equation}
Indeed, we have
\begin{align*}
\mu_{1}\left(  A-B\right)   &  =\mu_{1}\left(  \rho^{\prime}E_{n}\right)
=\rho^{\prime}n,\\
\mu_{k}\left(  A-B\right)   &  =\mu_{k}\left(  \rho^{\prime}E_{n}\right)
=0,\text{ for }k=2,...,n.
\end{align*}
Weyl's inequalities (e.g., see \cite{JoHo}, p. 181) imply, that if $C$ and $D$
are two Hermitian matrices of order $n$ then
\[
\mu_{2}\left(  C+D\right)  \leq\mu_{2}\left(  C\right)  +\mu_{1}\left(
D\right)  ,
\]
and
\[
\mu_{n}\left(  C+D\right)  \leq\mu_{n}\left(  C\right)  +\mu_{1}\left(
D\right)  .
\]
Hence, we see that
\[
\mu_{2}\left(  A\right)  \leq\mu_{1}\left(  B\right)  +\mu_{2}\left(
A-B\right)  =\mu_{1}\left(  B\right)  +\mu_{2}\left(  \rho E_{n}\right)
=\mu_{1}\left(  B\right)  ,
\]
and thus,
\begin{equation}
\mu_{2}\left(  A\right)  \leq\mu_{1}\left(  B\right)  \leq\sigma_{1}\left(
B\right)  . \label{eq2}%
\end{equation}
Similarly,
\[
\mu_{1}\left(  -B\right)  +\mu_{n}\left(  A\right)  \geq\mu_{n}\left(
A-B\right)  =\mu_{n}\left(  \rho E_{n}\right)  =0.
\]
and thus,
\[
\sigma_{1}\left(  B\right)  \geq\mu_{1}\left(  -B\right)  \geq-\mu_{n}\left(
A\right)  .
\]
This, together with (\ref{eq2}), implies (\ref{eq1}).

Let now $\mathbf{x}\in\mathbb{C}^{n}$ be a unit vector such that $\left\vert
\left\langle B\mathbf{x,x}\right\rangle \right\vert =\sigma_{1}\left(
B\right)  .$ Applying Lemma \ref{lapp} with $\varepsilon=1/3,$ we can find a
vector $\mathbf{y}=\left(  y_{i}\right)  _{1}^{n}\in\mathbb{C}^{n}$
satisfying
\[
\left\Vert \mathbf{x-y}\right\Vert \leq1/3,
\]
such that $y_{i}$ take $m$ distinct values $\alpha_{1}<...<\alpha_{m},$ where%
\[
m\leq\left\lceil \frac{8\pi}{1/3}\right\rceil \left\lceil \frac{4}{1/3}%
\log\frac{4n}{1/3}\right\rceil .
\]
For every $i\in\left[  m\right]  ,$ let%
\[
P_{i}=\left\{  j:y_{j}=\alpha_{i}\right\}  ;
\]
clearly, $\left[  n\right]  =P_{1}\cup...\cup P_{m}$ is a proper partition.

We shall prove that
\begin{equation}
\sigma_{1}\left(  B\right)  \leq\frac{9}{2}\left\vert \left\langle
B\mathbf{y,y}\right\rangle \right\vert . \label{eq4}%
\end{equation}
Indeed, we have
\begin{align*}
\left\langle B\mathbf{x,x}\right\rangle -\left\langle B\mathbf{y,y}%
\right\rangle  &  =\left\langle B\left(  \mathbf{y-x}\right)  ,\mathbf{x}%
\right\rangle +\overline{\left\langle \mathbf{y,B\left(  x\mathbf{-y}\right)
}\right\rangle }\\
&  \leq\left(  \left\Vert \mathbf{y}\right\Vert +\left\Vert \mathbf{x}%
\right\Vert \right)  \left\Vert B\left(  \mathbf{y-x}\right)  \right\Vert \\
&  \leq\left(  2+\frac{1}{3}\right)  \sigma_{1}\left(  B\right)  \left\Vert
\mathbf{x-y}\right\Vert \leq\frac{7}{9}\sigma_{1}\left(  B\right)  .
\end{align*}
Hence if $\left\langle B\mathbf{x,x}\right\rangle =\sigma_{1}\left(  B\right)
$ then (\ref{eq4}) holds. Also, we have%
\begin{align*}
\left\langle B\left(  \mathbf{y-x}\right)  ,\mathbf{x}\right\rangle
+\overline{\left\langle \mathbf{y,B\left(  x\mathbf{-y}\right)  }\right\rangle
}  &  \geq-\left(  \left\Vert \mathbf{y}\right\Vert +\left\Vert \mathbf{x}%
\right\Vert \right)  \left\Vert B\left(  \mathbf{y-x}\right)  \right\Vert
\geq\\
-\left(  2+\frac{1}{3}\right)  \sigma_{1}\left(  B\right)  \left\Vert
\mathbf{x-y}\right\Vert  &  \geq-\frac{7}{9}\sigma_{1}\left(  B\right)  ,
\end{align*}
and thus, (\ref{eq4}) holds also if $\left\langle B\mathbf{x,x}\right\rangle
=-\sigma_{1}\left(  B\right)  .$

Define the Hermitian matrix $C=\left(  c_{ij}\right)  _{i,j=1}^{m}$ by
\[
c_{ij}=\frac{1}{\sqrt{\left\vert P_{i}\right\vert \left\vert P_{j}\right\vert
}}\sum_{r\in P_{i}}\sum_{s\in P_{j}}b_{rs}.
\]
Applying Lemma \ref{le7} to the partition $\left[  n\right]  =P_{1}\cup...\cup
P_{m}$ and the vector $\mathbf{y}$ we find that
\[
\left\vert \left\langle B\mathbf{y,y}\right\rangle \right\vert \leq\sigma
_{1}\left(  C\right)  .
\]
Hence, in view of (\ref{eq1}) and (\ref{eq4}), we see that
\[
\sigma_{2}\left(  A\right)  \leq\sigma_{1}\left(  B\right)  \leq\frac{9}%
{2}\sigma_{1}\left(  C\right)  .
\]
Observe that,
\begin{align*}
\sigma_{1}\left(  C\right)   &  =\max\left(  \left\vert \mu_{1}\left(
C\right)  \right\vert ,\left\vert \mu_{n}\left(  C\right)  \right\vert
\right)  \leq m\max_{i,j\in\left[  m\right]  }\left\vert c_{ij}\right\vert \\
&  \leq\left\lceil \frac{8\pi}{1/3}\right\rceil \left\lceil \frac{4}{1/3}%
\log\frac{4n}{1/3}\right\rceil \max_{i,j\in\left[  m\right]  }\left\vert
c_{ij}\right\vert .
\end{align*}
Since,
\[
\max_{i,j\in\left[  m\right]  }\left\vert c_{ij}\right\vert \leq disc\left(
A\right)  ,
\]
we obtain%
\begin{equation}
\sigma_{2}\left(  A\right)  \leq\frac{9}{2}\left\lceil \frac{8\pi}%
{1/3}\right\rceil \left\lceil \frac{4}{1/3}\log\frac{4n}{1/3}\right\rceil
disc\left(  A\right)  , \label{exprC}%
\end{equation}
and the proof is completed.
\end{proof}

In the arguments above we made no attempt to optimize the constant in Theorem
\ref{th1}. As the right-hand side of (\ref{exprC}) is bounded above by
\[
\left(  4104\log n+10260\right)  disc\left(  A\right)  ,
\]
we can take $C$ to be $18906$.

\subsection{\label{sec2.2}Tightness of the upper bound on $\sigma_{2}$}

For $n=2k\geq2,$ let $A^{\prime}=\left(  a_{ij}^{\prime}\right)  _{i,j=1}^{k}$
be defined by%
\[
a_{ij}^{\prime}=\frac{1}{\sqrt{ij}},
\]
and let $A=\left(  a_{ij}\right)  _{i,j=1}^{n}$ be the block matrix
\[
A=\left(
\begin{array}
[c]{ll}%
E_{k}+A^{\prime} & E_{k}-A^{\prime}\\
E_{k}-A^{\prime} & E_{k}+A^{\prime}%
\end{array}
\right)  .
\]

Clearly, $A$ is nonnegative and symmetric. As we shall see the matrix $A$
shows that Theorem \ref{th1} is best possible up to a multiplicative constant.

\begin{theorem}
\label{th01} For the matrix $A$ defined above we have
\begin{equation}
\mu_{2}\left(  A\right)  \geq\frac{1}{2}disc\left(  A\right)  \log n.
\label{sigma2}%
\end{equation}

\end{theorem}

\begin{proof}
In fact, we shall show that $\mu_{2}\left(  A\right)  $ and $disc\left(
A\right)  $ satisfy
\[
\mu_{2}\left(  A\right)  \geq2\log n
\]
and%
\begin{equation}
disc\left(  A\right)  <4. \label{co2}%
\end{equation}
Indeed, the sum of every row of $A$ is exactly $n,$ and, since $A$ is
nonnegative, it follows that $\mu_{1}\left(  A\right)  =n.$ Note that the
vector $\mathbf{j}\in\mathbb{R}^{n}$ of all ones is an eigenvector of $A$ to
$\mu_{1}\left(  A\right)  .$ By the Rayleigh principle%
\[
\mu_{2}\left(  A\right)  =\max_{\mathbf{y\bot j,y\neq0}}\frac{\left\langle
A\mathbf{y,y}\right\rangle }{\left\Vert \mathbf{y}\right\Vert ^{2}},
\]
so our goal is to find a nonzero $\mathbf{y}\in\mathbb{R}^{n}$ such that
$\mathbf{y}\bot\mathbf{j,}$ and the ratio $\left\langle A\mathbf{y,y}%
\right\rangle /\left\Vert \mathbf{y}\right\Vert ^{2}$ is sufficiently large.

Define the vector $\mathbf{y}=\left(  y_{i}\right)  _{i=1}^{n}$ by
\[
y_{i}=\left\{
\begin{array}
[c]{lll}%
1/\sqrt{i} & \text{if} & i\leq k\\
-1/\sqrt{i-k} & \text{if} & i>k.
\end{array}
\right.  .
\]
From%
\[
\sum_{i=1}^{2k}y_{i}=\sum_{i=1}^{k}\frac{1}{\sqrt{i}}-\sum_{i=k+1}^{2k}%
\frac{1}{\sqrt{i-k}}=0
\]
we see that $\mathbf{y}\bot\mathbf{j.}$ Setting
\[
\xi_{k}=\sum_{i=1}^{k}\frac{1}{i},
\]
we deduce
\[
\left\Vert \mathbf{y}\right\Vert ^{2}=\sum_{i=1}^{k}\frac{1}{i}+\sum
_{i=k+1}^{n}\frac{1}{i-k}=2\sum_{i=1}^{k}\frac{1}{i}=2\xi_{k}.
\]
Next, we shall compute $\left\langle A\mathbf{y,y}\right\rangle .$ Recall
that
\[
a_{ij}=\left\{
\begin{array}
[c]{cccc}%
1+1/\sqrt{ij} & \text{if} & i\leq k, & j\leq k\\
1-1/\sqrt{ij} & \text{if} & i\leq k, & j>k\\
1-1/\sqrt{ij} & \text{if} & i>k, & j\leq k\\
1+1/\sqrt{ij} & \text{if} & i>k & j>k
\end{array}
\right.  .
\]
Thus we have
\begin{align*}
\left\langle A\mathbf{y,y}\right\rangle  &  =\sum_{i=1}^{2k}\sum_{j=1}%
^{2k}a_{ij}y_{i}y_{j}=\sum_{i=1}^{k}\sum_{j=1}^{k}\frac{a_{ij}}{\sqrt{ij}%
}+\sum_{i=k+1}^{2k}\sum_{j=k+1}^{2k}\frac{a_{ij}}{\sqrt{\left(  i-k\right)
\left(  j-k\right)  }}\\
&  -\sum_{i=1}^{k}\sum_{j=k+1}^{2k}\frac{a_{ij}}{\sqrt{i\left(  j-k\right)  }%
}-\sum_{i=k+1}^{2k}\sum_{j=1}^{k}\frac{a_{ij}}{\sqrt{\left(  i-k\right)  j}}\\
&  =\sum_{i=1}^{k}\sum_{j=1}^{k}\frac{1}{\sqrt{ij}}\left(  1+\frac{1}%
{\sqrt{ij}}\right)  +\sum_{i=1}^{k}\sum_{j=1}^{k}\frac{1}{\sqrt{ij}}\left(
1+\frac{1}{\sqrt{ij}}\right) \\
&  -\sum_{i=1}^{k}\sum_{j=1}^{k}\frac{1}{\sqrt{ij}}\left(  1-\frac{1}%
{\sqrt{ij}}\right)  -\sum_{i=1}^{k}\sum_{j=1}^{k}\frac{1}{\sqrt{ij}}\left(
1-\frac{1}{\sqrt{ij}}\right) \\
&  =4\sum_{i=1}^{k}\sum_{j=1}^{k}\frac{1}{\sqrt{ij}}\frac{1}{\sqrt{ij}%
}=4\left(  \xi_{k}\right)  ^{2}.
\end{align*}
Hence,%
\[
\mu_{2}\left(  A\right)  \geq\frac{\left\langle A\mathbf{y,y}\right\rangle
}{\left\Vert \mathbf{y}\right\Vert ^{2}}\geq2\xi_{k}>2\log n.
\]

Let us now turn to our proof of (\ref{co2}). Since the sum of every row of $A$
is exactly $n,$ we have $\rho^{\prime}\left(  A\right)  =1.$

Assume $X_{0},Y_{0}\subset\left[  n\right]  $ are nonempty sets, maximizing
the right-hand side of (\ref{defdisc}), i.e. satisfying
\begin{equation}
disc\left(  A\right)  =\frac{1}{\sqrt{\left\vert X_{0}\right\vert \left\vert
Y_{0}\right\vert }}\left\vert \sum_{i\in X_{0}}\sum_{j\in Y_{0}}\left(
a_{ij}-1\right)  \right\vert . \label{disceq}%
\end{equation}
Set
\begin{align*}
X_{1}  &  =X_{0}\cap\left[  k\right]  ,\text{ }X_{2}=X_{0}\cap\left[
k+1,n\right]  ,\\
Y_{1}  &  =Y_{0}\cap\left[  k\right]  ,\text{ }Y_{2}=Y_{0}\cap\left[
k+1,n\right]  .
\end{align*}
Then the right-hand side of (\ref{disceq}) is equal to%
\begin{align*}
&  \frac{1}{\sqrt{\left\vert X_{0}\right\vert \left\vert Y_{0}\right\vert }%
}\left\vert \sum_{i\in X_{1}}\sum_{j\in Y_{1}}\frac{1}{\sqrt{ij}}+\sum_{i\in
X_{2}}\sum_{j\in Y_{2}}\frac{1}{\sqrt{ij}}-\sum_{i\in X_{1}}\sum_{j\in Y_{2}%
}\frac{1}{\sqrt{ij}}-\sum_{i\in X_{2}}\sum_{j\in Y_{1}}\frac{1}{\sqrt{ij}%
}\right\vert \\
&  =\frac{1}{\sqrt{\left\vert X_{0}\right\vert \left\vert Y_{0}\right\vert }%
}\left\vert \left(  \sum_{i\in X_{1}}\frac{1}{\sqrt{i}}-\sum_{i\in X_{2}}%
\frac{1}{\sqrt{i}}\right)  \left(  \sum_{i\in Y_{1}}\frac{1}{\sqrt{i}}%
-\sum_{i\in Y_{2}}\frac{1}{\sqrt{i}}\right)  \right\vert
\end{align*}
Since $disc\left(  A\right)  $ is maximal, one of $X_{1},X_{2}$ is empty, and
one of $Y_{1},Y_{2}$ is empty. By symmetry we can assume that $X_{2}%
=\varnothing,$ $Y_{2}=\varnothing.$ Then the matrix $A\left[  X_{0}%
,Y_{0}\right]  =A\left[  X_{1},Y_{1}\right]  $ is in the upper-left-hand
corner of $A$ and
\begin{align*}
\frac{1}{\sqrt{\left\vert X_{0}\right\vert \left\vert Y_{0}\right\vert }%
}\left\vert \sum_{i\in X_{0}}\sum_{j\in Y_{0}}\left(  a_{ij}-1\right)
\right\vert  &  =\frac{1}{\sqrt{\left\vert X_{0}\right\vert \left\vert
Y_{0}\right\vert }}\left(  \sum_{i=1}^{\left\vert X_{0}\right\vert }\frac
{1}{\sqrt{i}}\right)  \left(  \sum_{i=1}^{\left\vert Y_{0}\right\vert }%
\frac{1}{\sqrt{i}}\right) \\
&  <\frac{4\sqrt{\left\vert X_{0}\right\vert }\sqrt{\left\vert Y_{0}%
\right\vert }}{\sqrt{\left\vert X_{0}\right\vert \left\vert Y_{0}\right\vert
}}=4.
\end{align*}

\end{proof}

It is not impossible that the constant 4 appearing in (\ref{co2}) is fairly
close to be the best possible.

\section{\label{sec3} A class of dense regular graphs}

Our goal in this section is to construct infinitely many regular graphs $G$
such that
\[
\mu_{2}\left(  G\right)  \geq Cdisc_{2}\left(  G\right)  \log\left(  \left|
G\right|  \right)
\]
for some absolute constant $C>0.$ In fact, for every sufficiently large prime
$p$ and $k=\left\lceil p^{1/5}\right\rceil $ we shall construct a matrix
$\mathcal{A}$ such that:

\emph{(a)} $\mathcal{A}$ is a square, symmetric, $\left(  0,1\right)  $-matrix
of size $2kp$ with zero main diagonal;

\emph{(b)} all row sums of $\mathcal{A}$ are equal to $kp;$

\emph{(c)} $\mu_{2}\left(  \mathcal{A}\right)  $ satisfies
\[
\mu_{2}\left(  \mathcal{A}\right)  \geq\frac{1}{2}p\log k;
\]

\emph{(d)} $disc\left(  \mathcal{A}\right)  $ satisfies
\[
disc\left(  \mathcal{A}\right)  \leq12p.
\]
The matrix $\mathcal{A}$ will be constructed as a block matrix of $4k^{2}$
blocks, each block being a square matrix of size $p.$

We shall select a symmetric matrix of integers that is roughly proportional to
the matrix $A$ of section \ref{sec2.2}, and then we shall replace each entry
of that matrix by a $p\times p,$ symmetric, $\left(  0,1\right)  $-matrix of
low discrepancy and density equal to the value of the corresponding entry.

Before describing the blocks of $\mathcal{A}$, we shall consider a corollary
of a theorem of Thomason.

\subsection{\label{Tho}A theorem of Thomason}

Thomason (\cite{Tho3}, Theorem 2) proved a widely-applicable result about
bipartite graphs with vertex classes of equal size; for convenience, we shall
restate his theorem in matrix form.

\begin{theorem}
Let $0<p<1,$ $\mu\geq0,$ and $A$ be a square $\left(  0,1\right)  $-matrix of
size $n$. If each row of $A$ has at least $pn$ ones, and the inner product of
every two distinct rows is at most $p^{2}n+\mu,$ then for every $X,Y\subset
\left[  n\right]  ,$%
\[
\left\vert \sum_{i\in X}\sum_{j\in Y}\left(  a_{ij}-p\right)  \right\vert
\leq\varepsilon\left\vert Y\right\vert +\sqrt{\left\vert X\right\vert
\left\vert Y\right\vert \left(  pn+\mu\left\vert X\right\vert \right)  },
\]
where $\varepsilon=1$ if $p\left\vert X\right\vert <1$ and $\varepsilon=0$ otherwise.
\end{theorem}

Applying this theorem to the adjacency matrix of a graph $G,$ we obtain
immediately the following generalization of Theorem 1 in \cite{Tho3}.

\begin{theorem}
\label{GenT}Let $0<p<1,$ $\mu\geq0,$ and $G$ be a graph of order $n.$ If
$d\left(  u\right)  \geq pn$ for every $u\in V\left(  G\right)  ,$ and
\[
\left\vert \Gamma\left(  u\right)  \cap\Gamma\left(  v\right)  \right\vert
\leq p^{2}n+\mu
\]
for every two distinct $u,v\in V\left(  G\right)  ,$ then for every
$X,Y\subset V\left(  G\right)  ,$%
\[
\left\vert e\left(  X,Y\right)  -p\left\vert X\right\vert \left\vert
Y\right\vert \right\vert \leq\varepsilon\left\vert Y\right\vert +\sqrt
{\left\vert X\right\vert \left\vert Y\right\vert \left(  pn+\mu\left\vert
X\right\vert \right)  },
\]
where $\varepsilon=1$ if $p\left\vert X\right\vert <1$ and $\varepsilon=0$ otherwise.
\end{theorem}

Next we shall describe a family of symmetric $\left(  0,1\right)  $-matrices
of size $p$ that we shall use as blocks of $\mathcal{A}.$

\subsection{The blocks of $\mathcal{A}$}

Let $p$ be a sufficiently large prime, $\mathbb{Z}_{p}$ be the field of order
$p,$ and $t\in\left[  p\right]  .$ Let $Q\left(  p,t\right)  $ be the graph
whose vertex set is $\left[  p\right]  ,$ and two distinct $u,v\in\left[
p\right]  $ are joined if%
\[
\left\{  \frac{\left(  u-v\right)  ^{2}}{p}\right\}  \leq\frac{t}{p},
\]
where $\left\{  x\right\}  $ is the fractional part of $x$. The graphs
$Q\left(  p,t\right)  $ were introduced by Bollob\'{a}s and Erd\H{o}s in
\cite{BE}, as examples of pseudo-random graphs. The following lemma summarizes
the properties of $Q\left(  p,t\right)  $ that we shall be interested in.

\begin{lemma}
\label{pqpt} The graph $Q\left(  p,t\right)  $ is a regular graph of order $p$
such that

\emph{(i) }the degree $d$ of $Q\left(  p,t\right)  $ satisfies%
\[
\left|  d-t\right|  \leq\sqrt{p}\left(  \log p\right)  ^{2};
\]

\emph{(ii)} the adjacency matrix $A$ of $Q\left(  p,t\right)  $ satisfies
\emph{ }%
\[
disc\left(  A\right)  <2p^{3/4}\log p.
\]

\end{lemma}

\begin{proof}
Since $Q\left(  p,t\right)  $ is invariant under the cyclic shift
$z\rightarrow z+1$ $\operatorname{mod}$ $p,$ it is clear that $Q\left(
p,t\right)  $ is regular. In fact, (i) follows from a much stronger result of
Burgess \cite{Bur}.

To prove (ii) we shall first recall that Theorem 3.16 in \cite{Bol2} states
that, for any two vertices $u,v$ of $Q\left(  p,t\right)  ,$ we have
\begin{equation}
\left\vert \left\vert \Gamma\left(  u\right)  \cap\Gamma\left(  v\right)
\right\vert -\frac{t^{2}}{p}\right\vert <\sqrt{p}\left(  \log p\right)  ^{2},
\label{nei1}%
\end{equation}
Setting $\beta=d/p,$ from (i) and (\ref{nei1}), for every two vertices $u,v$
of $Q\left(  p,t\right)  ,$ we obtain
\begin{align*}
\left\vert \Gamma\left(  u\right)  \cap\Gamma\left(  v\right)  \right\vert  &
\leq\frac{t^{2}}{p}+\sqrt{p}\left(  \log p\right)  ^{2}\leq\frac{\left(  \beta
p+\sqrt{p}\left(  \log p\right)  ^{2}\right)  ^{2}}{p}+\sqrt{p}\left(  \log
p\right)  ^{2}\\
&  =\beta^{2}p+2\beta\sqrt{p}\left(  \log p\right)  ^{2}+\left(  \log
p\right)  ^{4}+\sqrt{p}\left(  \log p\right)  ^{2}\\
&  <\beta^{2}p+3\sqrt{p}\left(  \log p\right)  ^{2}+\left(  \log p\right)
^{4}.
\end{align*}

Suppose that $X,Y\subset\left[  p\right]  $ are nonempty sets. Assuming
$\left\vert X\right\vert \leq\left\vert Y\right\vert ,$ by Theorem \ref{GenT},
we obtain%
\[
\left\vert e\left(  X,Y\right)  -\beta\left\vert X\right\vert \left\vert
Y\right\vert \right\vert \leq\left\vert X\right\vert +\sqrt{\left\vert
X\right\vert \left\vert Y\right\vert }\sqrt{\beta p+\left(  3\sqrt{p}\left(
\log p\right)  ^{2}+\left(  \log p\right)  ^{4}\right)  \left\vert
Y\right\vert }.
\]
Hence, noting that $\left\vert Y\right\vert \leq p$ and $\beta<1,$ we find
that%
\begin{align*}
\frac{1}{\sqrt{\left\vert X\right\vert \left\vert Y\right\vert }}\left\vert
e\left(  X,Y\right)  -\beta\left\vert X\right\vert \left\vert Y\right\vert
\right\vert  &  \leq1+\sqrt{\beta p+\left(  3\sqrt{p}\left(  \log p\right)
^{2}+\left(  \log p\right)  ^{4}\right)  p}\\
&  <2p^{3/4}\log p.
\end{align*}
Let $A=\left(  a_{ij}\right)  _{i,j=1}^{p}.$ Since, for every $X,Y\subset
\left[  p\right]  ,$ we have
\[
\sum_{i\in X}\sum_{j\in Y}a_{ij}=e\left(  X,Y\right)  ,
\]
and $\rho^{\prime}\left(  A\right)  =\beta,$ we deduce
\[
disc\left(  A\right)  <2p^{3/4}\log p,
\]
as claimed.
\end{proof}

Let $\mathcal{V}_{p}$ be the set of the degrees of the graphs $Q\left(
p,t\right)  $ for $t\in\left[  p\right]  .$ From Lemma \ref{pqpt}, (i), we see
that for every $s\in\left[  p\right]  $ there is a $d\in\mathcal{V}_{p},$ such
that there exists a $d$-regular graph $H\left(  p,d\right)  $ with%
\[
\left\vert d-s\right\vert \leq\sqrt{p}\log^{2}p,
\]
and
\[
disc_{2}\left(  H\left(  p,d\right)  \right)  <2p^{3/4}\log p.
\]

Now, for every $d\in\mathcal{V}_{p}$, let $A\left(  p,d\right)  $ be the
adjacency matrix of $H\left(  p,d\right)  .$ The properties of the matrices
$\left\{  A\left(  p,d\right)  :d\in\mathcal{V}_{p}\right\}  $ are summarized
in the following lemma.

\begin{lemma}
\label{paq} For every integer $s\in\left[  p\right]  ,$ there exist
$d\in\mathcal{V}_{p}$ and a matrix $A\left(  p,d\right)  ,$ such that

\emph{(i) }$\left|  d-s\right|  <\sqrt{p}\log^{2}p;$

\emph{(ii)} $A\left(  p,d\right)  $ is a symmetric $\left(  0,1\right)
$-matrix of size $p$ with zero main diagonal;

\emph{(iii) }all row sums of $A\left(  p,d\right)  $ are equal to $d;$

\emph{(iv)} the function $disc\left(  A\left(  p,d\right)  \right)  $
satisfies
\[
disc\left(  A\left(  p,d\right)  \right)  <2p^{3/4}\log p.
\]

\end{lemma}

If $A$ is a square $\left(  0,1\right)  $-matrix of size $n,$ we call the
matrix
\[
\overline{A}=E_{n}-A
\]
the \emph{complement} of $A$. Observe that if $A$ is a square $\left(
0,1\right)  $-matrix then%
\begin{align*}
\rho^{\prime}\left(  \overline{A}\right)   &  =1-\rho^{\prime}\left(
A\right)  ,\\
disc\left(  \overline{A}\right)   &  =disc\left(  A\right)  .
\end{align*}
Hence, the complement of any matrix $A\left(  p,d\right)  $ satisfies
\[
disc\left(  \overline{A\left(  p,d\right)  }\right)  <2p^{3/4}\log p.
\]
The matrices $\left\{  A\left(  p,d\right)  :d\in\mathcal{V}_{p}\right\}  $
together with their complements will be used as blocks of the matrix
$\mathcal{A}.$

\subsection{The construction of $\mathcal{A}$}

For every $s\in\left[  2k\right]  ,$ set
\[
I_{s}=\left\{  i:\left(  s-1\right)  p\leq i<sp\right\}  .
\]
Define the matrix $D=\left(  d_{ij}\right)  _{i,j=1}^{k}$ by%
\[
d_{ij}=q,\text{ }q\in\mathcal{V}_{p}\text{, }\left\vert q-\left(  \frac{p}%
{2}+\frac{p}{2\sqrt{ij}}\right)  \right\vert =\min_{x\in\mathcal{V}_{p}%
}\left\vert x-\left(  \frac{p}{2}+\frac{p}{2\sqrt{ij}}\right)  \right\vert .
\]
From Lemma \ref{pqpt}, (i), we see that%
\begin{equation}
\left\vert 2d_{ij}-\left(  p+\frac{p}{\sqrt{ij}}\right)  \right\vert
\leq2\sqrt{p}\log^{2}p. \label{dij}%
\end{equation}
The matrix $D$ will be the cornerstone of our construction. Note that $D$ is
symmetric and the values of its entries belong to the set $\mathcal{V}_{p}.$

Now, let us define $\mathcal{A}^{\prime}$ as a block matrix by
\begin{equation}
\mathcal{A}^{\prime}=\left(
\begin{array}
[c]{cccc}%
A\left(  p,d_{11}\right)  & A\left(  p,d_{12}\right)  & . & A\left(
p,d_{1k}\right) \\
A\left(  p,d_{12}\right)  & A\left(  p,d_{22}\right)  & . & .\\
. & . & . & .\\
A\left(  p,d_{1k}\right)  & . & . & A\left(  p,d_{kk}\right)
\end{array}
\right)  , \label{defA1}%
\end{equation}
and set
\begin{equation}
\mathcal{A}=\left(
\begin{array}
[c]{cc}%
\mathcal{A}^{\prime} & E_{kp}-\mathcal{A}^{\prime}\\
E_{kp}-\mathcal{A}^{\prime} & \mathcal{A}^{\prime}%
\end{array}
\right)  . \label{defA}%
\end{equation}
By our construction $\mathcal{A}$ is a symmetric $\left(  0,1\right)  $-matrix
of size $2pk,$ and its main diagonal is zero, so $\mathcal{A}$ satisfies (a).
Also, we see that every row sum of $\mathcal{A}$ is exactly $kp$, so
$\mathcal{A}$ satisfies (b) as well. In the following two theorems we shall
prove that $\mathcal{A}$ satisfies also (c) and (d).

For the sake of convenience, set $\mathcal{A}=\left(  a_{ij}\right)
_{i,i=1}^{2pk}$ and $\mathcal{A}_{ij}=\mathcal{A}\left[  I_{i},I_{j}\right]  $
for $i,j\in\left[  2k\right]  .$ Observe that the row sums of any matrix
$\mathcal{A}_{ij}$ are equal, and from Lemma \ref{paq} and what follows, we
have
\begin{equation}
disc\left(  \mathcal{A}_{ij}\right)  \leq2p^{3/4}\log p. \label{daij}%
\end{equation}

\begin{theorem}
\label{prop2} The second eigenvalue $\mu_{2}\left(  \mathcal{A}\right)  $ of
the matrix $\mathcal{A}$ defined by (\ref{defA}) satisfies
\[
\mu_{2}\left(  \mathcal{A}\right)  \geq\frac{1}{2}p\log k.
\]

\end{theorem}

\begin{proof}
Indeed, from (\ref{defA}) we see that the sum of every row of $\mathcal{A}$ is
exactly $kp$. Since $\mathcal{A}$ is nonnegative, it follows that $\mu
_{1}\left(  \mathcal{A}\right)  =pk$ and the vector $\mathbf{j}\in
\mathbb{R}^{2pk}$ of all ones is an eigenvector of $\mathcal{A}$ to $\mu
_{1}\left(  \mathcal{A}\right)  .$ By the Rayleigh principle%
\[
\mu_{2}\left(  \mathcal{A}\right)  =\max_{\mathbf{y}\bot\mathbf{j,y\neq0}%
}\frac{\left\langle A\mathbf{y,y}\right\rangle }{\left\Vert \mathbf{y}%
\right\Vert ^{2}},
\]
so our goal is to find a nonzero $\mathbf{y}\in\mathbb{R}^{2pk}$ such that
$\mathbf{y}\bot\mathbf{j}$ and the ratio $\left\langle A\mathbf{y,y}%
\right\rangle /\left\Vert \mathbf{y}\right\Vert ^{2}$ is sufficiently large.

Define the vector $\mathbf{y}=\left(  y_{i}\right)  _{i=1}^{n}$ by
\[
y_{i}=\left\{
\begin{array}
[c]{lll}%
1/\sqrt{s} & \text{if} & i\in I_{s},\text{ }s\leq k\\
-1/\sqrt{s-k} & \text{if} & i\in I_{s},\text{ }s>k.
\end{array}
\right.  .
\]
From%
\[
\sum_{i=1}^{2pk}y_{i}=\sum_{s=1}^{2k}\sum_{i\in I_{s}}\frac{1}{\sqrt{s}}%
-\sum_{s=k+1}^{2k}\sum_{i\in I_{s}}\frac{1}{\sqrt{s-k}}=\sum_{s=1}^{k}\frac
{p}{\sqrt{s}}-\sum_{s=1}^{k}\frac{p}{\sqrt{s}}=0
\]
we see that $\mathbf{y}\bot\mathbf{j.}$ Also, for $\left\Vert \mathbf{y}%
\right\Vert ^{2}$ we have
\[
\left\Vert \mathbf{y}\right\Vert ^{2}=\sum_{s=1}^{k}\sum_{i\in I_{s}}\frac
{1}{s}+\sum_{s=k+1}^{2k}\sum_{i\in I_{s}}\frac{1}{s-k}=2\sum_{s=1}^{k}%
\sum_{i\in I_{s}}\frac{1}{s}=2p\sum_{s=1}^{k}\frac{1}{s}=2p\xi_{k}.
\]
On the other hand, for $\left\langle \mathcal{A}\mathbf{y,y}\right\rangle $ we
see that
\[
\left\langle \mathcal{A}\mathbf{y,y}\right\rangle =\sum_{i=1}^{2pk}\sum
_{j=1}^{2pk}a_{ij}y_{i}y_{j}=\sum_{i=1}^{2k}\sum_{j=1}^{2k}\sum_{s\in I_{i}%
}\sum_{t\in I_{j}}a_{st}y_{s}y_{t}.
\]
By (\ref{defA1}) and (\ref{defA}), we have
\[
\sum_{s\in I_{i}}\sum_{t\in I_{j}}a_{st}=\left\{
\begin{array}
[c]{llll}%
pd_{ij} & \text{if} & i\leq k & j\leq k\\
p\left(  p-d_{i\left(  j-k\right)  }\right)  & \text{if} & i\leq k & j>k\\
p\left(  p-d_{\left(  i-k\right)  j}\right)  & \text{if} & i>k & j\leq k\\
pd_{\left(  i-k\right)  \left(  j-k\right)  } & \text{if} & i>k & j>k
\end{array}
\right.  .
\]
Hence, by the choice of $\mathbf{y}$,%
\begin{align*}
\left\langle \mathcal{A}\mathbf{y,y}\right\rangle  &  =\sum_{i=1}^{k}%
\sum_{j=1}^{k}\frac{pd_{ij}}{\sqrt{ij}}+\sum_{i=k+1}^{2k}\sum_{j=k+1}%
^{2k}\frac{pd_{\left(  i-k\right)  \left(  j-k\right)  }}{\sqrt{\left(
i-k\right)  \left(  j-k\right)  }}\\
&  -\sum_{i=1}^{k}\sum_{j=k+1}^{2k}\frac{p\left(  p-d_{i\left(  j-k\right)
}\right)  }{\sqrt{i\left(  j-k\right)  }}-\sum_{i=k+1}^{2k}\sum_{j=1}^{k}%
\frac{p\left(  p-d_{\left(  i-k\right)  j}\right)  }{\sqrt{\left(  i-k\right)
j}}\\
&  =2p\left(  \sum_{i=1}^{k}\sum_{j=1}^{k}\frac{d_{ij}}{\sqrt{ij}}-\sum
_{i=1}^{k}\sum_{j=1}^{k}\frac{p-d_{ij}}{\sqrt{ij}}\right)  =2p\left(
\sum_{i=1}^{k}\sum_{j=1}^{k}\frac{2d_{ij}-p}{\sqrt{ij}}\right)  .
\end{align*}

From (\ref{dij}), we have%
\[
\frac{2d_{ij}-p}{\sqrt{ij}}>\frac{1}{\sqrt{ij}}\left(  \frac{p}{\sqrt{ij}%
}-2\sqrt{p}\left(  \log p\right)  ^{2}\right)  =\frac{p}{ij}-\frac{2\sqrt
{p}\left(  \log p\right)  ^{2}}{\sqrt{ij}},
\]
and so,
\begin{align*}
\left\langle \mathcal{A}\mathbf{y,y}\right\rangle  &  >2p^{2}\left(
\sum_{i=1}^{k}\sum_{j=1}^{k}\frac{1}{ij}\right)  -4p\sqrt{p}\left(  \log
p\right)  ^{2}\sum_{i=1}^{k}\sum_{j=1}^{k}\frac{1}{\sqrt{ij}}\\
&  >2p^{2}\left(  \xi_{k}\right)  ^{2}-16kp\sqrt{p}\left(  \log p\right)
^{2}.
\end{align*}
Hence, as $k\leq p^{1/5}$ and $p$ is large, $\left\langle \mathcal{A}%
\mathbf{y,y}\right\rangle >p^{2}\left(  \xi_{k}\right)  ^{2},$ and thus,
\[
\mu_{2}\left(  \mathcal{A}\right)  \geq\frac{\left\langle \mathcal{A}%
\mathbf{y,y}\right\rangle }{\left\|  \mathbf{y}\right\|  ^{2}}\geq\frac{1}%
{2}p\xi_{k}>\frac{1}{2}p\log k
\]
as claimed.
\end{proof}

\begin{theorem}
\label{prop3} If $p$ is large, $disc\left(  \mathcal{A}\right)  $ of the
matrix $\mathcal{A}$ defined by (\ref{defA}) satisfies%
\[
disc\left(  \mathcal{A}\right)  \leq12p.
\]

\end{theorem}

\begin{proof}
Since all row sums of $\mathcal{A}$ are exactly $pk,$ we deduce $\rho^{\prime
}\left(  \mathcal{A}\right)  =1/2.$

As before, assume $X_{0},Y_{0}\subset\left[  2kp\right]  $ are nonempty sets,
maximizing the right-hand side of (\ref{defdisc}), i.e. satisfying
\begin{equation}
disc\left(  \mathcal{A}\right)  =\frac{1}{\sqrt{\left\vert X_{0}\right\vert
\left\vert Y_{0}\right\vert }}\left\vert \sum_{i\in X_{0}}\sum_{j\in Y_{0}%
}\left(  a_{ij}-\frac{1}{2}\right)  \right\vert . \label{disceq1}%
\end{equation}
Set
\[
J_{1}=\left[  kp\right]  ,\text{ }J_{2}=\left[  kp+1,2kp\right]
\]
and let
\[
X_{i}=X_{0}\cap J_{i},\text{ }Y_{i}=Y_{0}\cap J_{i},\text{ }i=1,2.
\]
For $i,j=1,2$ consider the value%
\[
\Delta_{ij}=\max_{X\subset J_{i},Y\subset J_{j},X\neq\varnothing
,Y\neq\varnothing}\frac{1}{\sqrt{\left\vert X\right\vert \left\vert
Y\right\vert }}\left\vert \sum_{i\in X}\sum_{j\in Y}\left(  a_{ij}-\frac{1}%
{2}\right)  \right\vert
\]
By (\ref{defA}), we have
\begin{align*}
\mathcal{A}\left[  J_{1},J_{1}\right]   &  =\mathcal{A}\left[  J_{2}%
,J_{2}\right]  ,\\
\mathcal{A}\left[  J_{1},J_{2}\right]   &  =\mathcal{A}\left[  J_{2}%
,J_{1}\right]  =E_{kn}-\mathcal{A}\left[  J_{1},J_{1}\right]  ,
\end{align*}
and hence,
\[
\Delta_{11}=\Delta_{12}=\Delta_{21}=\Delta_{22}.
\]
Consequently%
\begin{align*}
disc\left(  \mathcal{A}\right)   &  =\frac{1}{\sqrt{\left\vert X_{0}%
\right\vert \left\vert Y_{0}\right\vert }}\left\vert \sum_{i=1}^{2}\sum
_{j=1}^{2}\sum_{s\in X_{i}}\sum_{t\in Y_{j}}\left(  a_{ij}-\frac{1}{2}\right)
\right\vert \\
&  \leq\frac{1}{\sqrt{\left\vert X_{0}\right\vert \left\vert Y_{0}\right\vert
}}\sum_{i=1}^{2}\sum_{j=1}^{2}\Delta_{ij}\sqrt{\left\vert X_{i}\right\vert
\left\vert Y_{j}\right\vert }\\
&  =\frac{\Delta_{11}}{\sqrt{\left\vert X_{0}\right\vert \left\vert
Y_{0}\right\vert }}\left(  \sqrt{\left\vert X_{1}\right\vert }+\sqrt
{\left\vert X_{2}\right\vert }\right)  \left(  \sqrt{\left\vert Y_{1}%
\right\vert }+\sqrt{\left\vert Y_{2}\right\vert }\right) \\
&  \leq2\Delta_{11}\frac{\sqrt{\left(  \left\vert X_{1}\right\vert +\left\vert
X_{2}\right\vert \right)  \left(  \left\vert Y_{1}\right\vert +\left\vert
Y_{2}\right\vert \right)  }}{\sqrt{\left(  \left\vert X_{1}\right\vert
+\left\vert X_{2}\right\vert \right)  \left(  \left\vert Y_{1}\right\vert
+\left\vert Y_{2}\right\vert \right)  }}=2\Delta_{11}.
\end{align*}

To complete our proof we shall show that
\[
\Delta_{11}<6p.
\]
Fix some nonempty sets $X_{0},Y_{0}\subset\left[  kp\right]  $ such that%
\begin{equation}
\Delta_{11}=\frac{1}{\sqrt{\left\vert X_{0}\right\vert \left\vert
Y_{0}\right\vert }}\left\vert \sum_{i=1}^{k}\sum_{j=1}^{k}\sum_{s\in X_{i}%
}\sum_{t\in Y_{j}}\left(  a_{st}-\frac{1}{2}\right)  \right\vert ,
\label{delta11}%
\end{equation}
and for every $i\in\left[  k\right]  ,$ set%
\[
X_{i}=X_{0}\cap I_{i},\text{ }Y_{i}=Y_{0}\cap I_{i}.
\]
Observe that for $i,j\in\left[  k\right]  $ we have%
\[
\rho^{\prime}\left(  \mathcal{A}_{ij}\right)  -\frac{1}{2}=\frac{d_{ij}}%
{p}-\frac{1}{2},
\]
hence, by (\ref{dij}),
\[
\left\vert \rho^{\prime}\left(  \mathcal{A}_{ij}\right)  -\frac{1}%
{2}\right\vert <\frac{1}{2\sqrt{ij}}+\frac{2\left(  \log p\right)  ^{2}}%
{\sqrt{p}},
\]
and so,
\begin{align*}
\left\vert \sum_{s\in X_{i}}\sum_{t\in Y_{j}}\left(  a_{st}-\frac{1}%
{2}\right)  \right\vert  &  \leq\left\vert \sum_{s\in X_{i}}\sum_{t\in Y_{j}%
}\left(  a_{st}-\rho^{\prime}\left(  \mathcal{A}_{ij}\right)  \right)
\right\vert +\left\vert X_{i}\right\vert \left\vert Y_{j}\right\vert
\left\vert \rho^{\prime}\left(  \mathcal{A}_{ij}\right)  -\frac{1}%
{2}\right\vert \\
&  \leq disc\left(  \mathcal{A}_{ij}\right)  \sqrt{\left\vert X_{i}\right\vert
\left\vert Y_{j}\right\vert }+\left\vert X_{i}\right\vert \left\vert
Y_{j}\right\vert \left(  \frac{1}{2\sqrt{ij}}+\frac{2\left(  \log p\right)
^{2}}{\sqrt{p}}\right)  .
\end{align*}
Recalling (\ref{delta11}), we see that%
\begin{align}
\Delta_{11}  &  \leq\frac{1}{\sqrt{\left\vert X_{0}\right\vert \left\vert
Y_{0}\right\vert }}\sum_{i=1}^{k}\sum_{j=1}^{k}\left\vert \sum_{s\in X_{i}%
}\sum_{t\in Y_{j}}\left(  a_{st}-\frac{1}{2}\right)  \right\vert \nonumber\\
&  \leq\frac{1}{\sqrt{\left\vert X_{0}\right\vert \left\vert Y_{0}\right\vert
}}\sum_{i=1}^{k}\sum_{j=1}^{k}disc\left(  \mathcal{A}_{ij}\right)
\sqrt{\left\vert X_{i}\right\vert \left\vert Y_{j}\right\vert }\label{t1}\\
&  +\frac{1}{2\sqrt{\left\vert X_{0}\right\vert \left\vert Y_{0}\right\vert }%
}\sum_{i=1}^{k}\sum_{j=1}^{k}\frac{\left\vert X_{i}\right\vert \left\vert
Y_{j}\right\vert }{\sqrt{ij}}\label{t2}\\
&  +\left(  \frac{2\left(  \log p\right)  ^{2}}{\sqrt{p}}\right)  \frac
{1}{\sqrt{\left\vert X_{0}\right\vert \left\vert Y_{0}\right\vert }}\sum
_{i=1}^{k}\sum_{j=1}^{k}\left\vert X_{i}\right\vert \left\vert Y_{j}%
\right\vert \label{t3}\\
&  =A+B+C.\nonumber
\end{align}
We shall estimate the terms (\ref{t1}), (\ref{t2}) and (\ref{t3}) separately.

From (\ref{daij}) we obtain%
\[
A\leq\frac{2p^{3/4}\log p}{\sqrt{\left\vert X_{0}\right\vert \left\vert
Y_{0}\right\vert }}\sum_{i=1}^{k}\sum_{j=1}^{k}\sqrt{\left\vert X_{i}%
\right\vert \left\vert Y_{j}\right\vert }.
\]
Hence, by%
\begin{align*}
\frac{1}{\sqrt{\left\vert X_{0}\right\vert \left\vert Y_{0}\right\vert }}%
\sum_{i=1}^{k}\sum_{j=1}^{k}\sqrt{\left\vert X_{i}\right\vert \left\vert
Y_{j}\right\vert }  &  =\frac{1}{\sqrt{\left\vert X_{0}\right\vert }}\left(
\sum_{i=1}^{k}\sqrt{\left\vert X_{i}\right\vert }\right)  \frac{1}%
{\sqrt{\left\vert Y_{0}\right\vert }}\left(  \sum_{j=1}^{k}\sqrt{\left\vert
Y_{j}\right\vert }\right) \\
&  \leq\sqrt{k}\sqrt{k}=k,
\end{align*}
we have%
\begin{equation}
A\leq2kp^{3/4}\log p\leq p. \label{bt1}%
\end{equation}
Next we turn to (\ref{t2}). Obviously,
\begin{equation}
B=\frac{1}{2\sqrt{\left\vert X_{0}\right\vert }}\left(  \sum_{i=1}^{k}%
\frac{\left\vert X_{i}\right\vert }{\sqrt{i}}\right)  \frac{1}{\sqrt
{\left\vert Y_{0}\right\vert }}\left(  \sum_{i=1}^{k}\frac{\left\vert
Y_{i}\right\vert }{\sqrt{i}}\right)  . \label{in10}%
\end{equation}
We shall show that
\begin{equation}
\frac{1}{\sqrt{\left\vert X_{0}\right\vert }}\left(  \sum_{i=1}^{k}%
\frac{\left\vert X_{i}\right\vert }{\sqrt{i}}\right)  \leq2\sqrt{2p}.
\label{xs}%
\end{equation}
Indeed, set
\[
s=\left\lfloor \frac{\left\vert X_{0}\right\vert }{p}\right\rfloor ,
\]
and observe that the left-hand side of (\ref{xs}) attains its maximum when
\begin{align*}
\left\vert X_{i}\right\vert  &  =p,\text{ \ \ \ }1\leq i\leq s,\\
\left\vert X_{s+1}\right\vert  &  =\left\vert X_{0}\right\vert -ps,\\
\left\vert X_{i}\right\vert  &  =0,\text{ \ \ \ }s+1<i\leq2k.
\end{align*}
Obviously (\ref{xs}) holds if $s=0,$ so we shall assume $s\geq1.$ Then we
have,
\[
\frac{1}{\sqrt{\left\vert X_{0}\right\vert }}\left(  \sum_{i=1}^{k}%
\frac{\left\vert X_{i}\right\vert }{\sqrt{i}}\right)  \leq\frac{1}%
{\sqrt{\left\vert X_{0}\right\vert }}\sum_{i=1}^{s+1}\frac{p}{\sqrt{i}}%
\leq\frac{2p\sqrt{s+1}}{\sqrt{\left\vert X_{0}\right\vert }}\leq\sqrt
{2p\frac{sp}{\left\vert X_{0}\right\vert }}\leq2\sqrt{2p}%
\]
and (\ref{xs}) follows.

Similarly, we see that
\[
\frac{1}{\sqrt{\left\vert Y_{0}\right\vert }}\left(  \sum_{i=1}^{k}%
\frac{\left\vert Y_{i}\right\vert }{\sqrt{i}}\right)  \leq2\sqrt{2p}%
\]
and hence, in view of (\ref{in10}), we find
\begin{equation}
B\leq4p. \label{bt2}%
\end{equation}
Finally,
\begin{equation}
C=\left(  \frac{2\left(  \log p\right)  ^{2}}{\sqrt{p}}\right)  \frac
{\left\vert X\right\vert \left\vert Y\right\vert }{\sqrt{\left\vert
X\right\vert \left\vert Y\right\vert }}\leq\left(  \frac{2\left(  \log
p\right)  ^{2}}{\sqrt{p}}\right)  \sqrt{kp}<p. \label{bt3}%
\end{equation}
Now, replacing (\ref{t1}), (\ref{t2}), (\ref{t3}) by (\ref{bt1}), (\ref{bt2}),
(\ref{bt3}), we obtain%
\[
\Delta_{11}<6p,
\]
and the proof is completed.
\end{proof}

\subsection{\label{sec3.last} A conjecture of Chung}

In \cite{Ch2} Chung studies a version of the Laplacian matrix a graph $G$ that
she denotes by $\mathcal{L}\left(  G\right)  .$ If $G$ is $d$-regular of order
$n$ the matrix $\mathcal{L}\left(  G\right)  $ is given by%
\begin{equation}
\mathcal{L}\left(  G\right)  =I_{n}-\frac{1}{d}A, \label{defL}%
\end{equation}
where $A$ is the adjacency matrix of $G\left(  n\right)  .$ Following Chung's
notation, the eigenvalues of $\mathcal{L}\left(  G\right)  $ are $\lambda
_{0}\leq...\leq\lambda_{n-1},$ with $\lambda_{0}=0.$

Set $\overline{\lambda}=\max_{i\neq0}\left|  1-\lambda_{i}\right|  ,$ and for
every $X\subset V\left(  G\right)  ,$ let $vol$ $X=\sum_{v\in X}d\left(
v\right)  .$ Chung asked the following question.

Let $G$ be a nonempty graph and $\alpha>0$ is such that if $X,Y\subset
V=V\left(  G\right)  $ then
\begin{equation}
\left\vert e\left(  X,Y\right)  -\frac{vol\text{ }X\text{ }vol\text{ }%
Y}{vol\text{ }V}\right\vert \leq\alpha\frac{\sqrt{vol\text{ }X\text{
}vol\text{ }Y\text{ }vol\text{ }\left(  V\backslash X\right)  \text{
}vol\text{ }\left(  V\backslash Y\right)  }}{vol\text{ }V}. \label{vols}%
\end{equation}
Is there an absolute constant $C$ such that $\overline{\lambda}\leq C\alpha$?

We shall check that the graph $G_{p},$ whose adjacency matrix $\mathcal{A}%
_{p}=\mathcal{A}$ we constructed in the previous section, answers this
question in the negative. Indeed, recall that $G_{p}$ is $kp$-regular graph of
order $n=2kp.$ Theorem \ref{prop3} implies that%
\[
\left|  e\left(  X,Y\right)  -\frac{\left|  X\right|  \left|  Y\right|  }%
{2}\right|  \leq\frac{C}{k}\sqrt{\left|  X\right|  \left|  Y\right|  \left(
n-\left|  X\right|  \right)  \left(  n-\left|  Y\right|  \right)  }%
\]
for some absolute constant $C>0,$ so (\ref{vols}) holds with $\alpha=C/k.$ By
(\ref{defL}) and Theorem \ref{prop2}, we see that
\[
\overline{\lambda}\geq\left|  1-\lambda_{1}\right|  \geq1-\left(  1-\frac
{\mu_{2}}{\mu_{1}}\right)  =\frac{\mu_{2}}{\mu_{1}}\geq\frac{p\log k}%
{2kp}=\frac{\log k}{2k}%
\]
and $\overline{\lambda}$ is greater than any fixed multiple of $\alpha$.

\section{\label{sec4} Sparse graphs with low discrepancy and high second
eigenvalue}

In \cite{CG3} Chung and Graham extend quasi-random properties to sparse
graphs, i.e., graphs $G\left(  n,m\right)  $ with $m=o\left(  n^{2}\right)  $.
Their approach is based on the following. Fix a function $p=p\left(  n\right)
$ with $0<p<1$ and
\[
\lim_{n\rightarrow\infty}pn=\infty.
\]
Let $\mathcal{G}_{p}$ be an infinite family of graphs $\left\{  G\left(
n\right)  :n\rightarrow\infty\right\}  $ such that, for every $G\left(
n\right)  \in\mathcal{G}_{p},$%
\begin{equation}
e\left(  G\left(  n\right)  \right)  =\left(  1+o\left(  1\right)  \right)
p\binom{n}{2}. \label{defGp}%
\end{equation}
Chung and Graham investigated a number of properties that a family
$\mathcal{G}_{p}$ can have; we shall be concerned with the following two here
(\cite{CG3}, p. 220):

\textbf{DISC}(1): For every $G\left(  n\right)  \in\mathcal{G}_{p},$ and for
all $X,Y\subset V\left(  G\right)  ,$%
\[
\left|  e\left(  X,Y\right)  -p\left|  X\right|  \left|  Y\right|  \right|
=o\left(  pn^{2}\right)  .
\]

\textbf{EIG}: For every $G\left(  n\right)  \in\mathcal{G}_{p},$
\[
\mu_{1}\left(  G\right)  =\left(  1+o\left(  1\right)  \right)  pn,\text{ and
}\sigma_{2}\left(  G\right)  =o\left(  pn\right)  .
\]

Chung and Graham proved that \textbf{EIG} implies \textbf{DISC}(1) (Theorem 1
in \cite{CG3}), and asked the following natural question (\cite{CG3}, p. 230).

\textbf{Question }Does \textbf{DISC}(1) imply \textbf{EIG}?

Recently Krivelevich and Sudakov (\cite{KrSu}, p. 9,) constructed an example
that answers this question in the negative. To conclude the paper we give a
general construction that we believe sheds more light on the relationship
between \textbf{DISC}(1) and \textbf{EIG}.

\begin{proposition}
For $p=p\left(  n\right)  =o\left(  1\right)  $ let $\mathcal{G}_{p}$ be a
family of graphs having the property \textbf{EIG}. Let $\mathcal{G}_{p}^{\ast
}$ be the family of the graphs that can be represented as disjoint unions
\[
G\left(  n\right)  \cup K_{\left\lfloor pn\right\rfloor },
\]
where $G\left(  n\right)  \in\mathcal{G}_{p}.$ Then $\mathcal{G}_{p}^{\ast}$
has \textbf{DISC}(1) but does not have \textbf{EIG.}.
\end{proposition}

\begin{proof}
Note that
\[
e\left(  G\left(  n\right)  \cup K_{\left\lfloor pn\right\rfloor }\right)
=\left(  1+o\left(  1\right)  \right)  p\binom{n}{2}+\binom{\left\lfloor
pn\right\rfloor }{2}=\left(  1+o\left(  1\right)  \right)  p\binom
{n+\left\lfloor pn\right\rfloor }{2},
\]
so $\mathcal{G}_{p}^{\ast}$ is defined according to (\ref{defGp}). Also, given
$G^{\prime}=G\left(  n\right)  \cup K_{\left\lfloor pn\right\rfloor },$
$Z=V\left(  K_{\left\lfloor pn\right\rfloor }\right)  $ and $X,Y\subset
V\left(  G^{\prime}\right)  ,$ we have
\begin{align*}
\left\vert e\left(  X,Y\right)  -p\left\vert X\right\vert \left\vert
Y\right\vert \right\vert  &  \leq\left\vert e\left(  X\backslash Z,Y\backslash
Z\right)  -p\left\vert X\backslash Z\right\vert \left\vert Y\backslash
Z\right\vert \right\vert \\
&  +\left\vert e\left(  X,Y\right)  -e\left(  X\backslash Z,Y\backslash
Z\right)  \right\vert +\left\vert p\left\vert X\right\vert \left\vert
Y\right\vert -p\left\vert X\backslash Z\right\vert \left\vert Y\backslash
Z\right\vert \right\vert \\
&  \leq o\left(  pn^{2}\right)  +2e\left(  Z\right)  +p\left\vert Z\right\vert
\left(  \left\vert X\right\vert +\left\vert Y\right\vert \right) \\
&  \leq o\left(  pn^{2}\right)  +p^{2}n^{2}+2p^{2}n^{2}=o\left(
pn^{2}\right)  .
\end{align*}
Thus, $\mathcal{G}_{p}^{\ast}$ has \textbf{DISC}(1). However, since
$G^{\prime}$ is a union of the disjoint graphs $G\left(  n\right)  $ and
$K_{\left\lfloor pn\right\rfloor },$ we find that
\begin{align*}
\min\left\{  \mu_{1}\left(  G\left(  n\right)  \right)  ,\mu_{1}\left(
K_{\left\lfloor pn\right\rfloor }\right)  \right\}   &  \leq\mu_{2}\left(
G^{\prime}\right)  \leq\mu_{1}\left(  G^{\prime}\left(  n\right)  \right) \\
&  =\max\left\{  \mu_{1}\left(  G\left(  n\right)  \right)  ,\mu_{1}\left(
K_{\left\lfloor pn\right\rfloor }\right)  \right\}  .
\end{align*}
Hence, from $\mu_{1}\left(  G\left(  n\right)  \right)  =\left(  1+o\left(
1\right)  \right)  pn$ and $\mu_{1}\left(  K_{\left\lfloor pn\right\rfloor
}\right)  =\left\lfloor pn\right\rfloor -1,$ we see that
\[
\mu_{2}\left(  G^{\prime}\right)  =\left(  1+o\left(  1\right)  \right)  pn,
\]
and so, $\mathcal{G}_{p}^{\ast}$ does not have \textbf{EIG}.
\end{proof}

\textbf{Acknowledgement} We are grateful to the referees for their valuable comments.

\end{document}